# 1st International Conference on Engineering and Applied Sciences Optimization

## PROCEEDINGS

M.G. Karlaftis, N.D. Lagaros, M. Papadrakakis (Eds.)

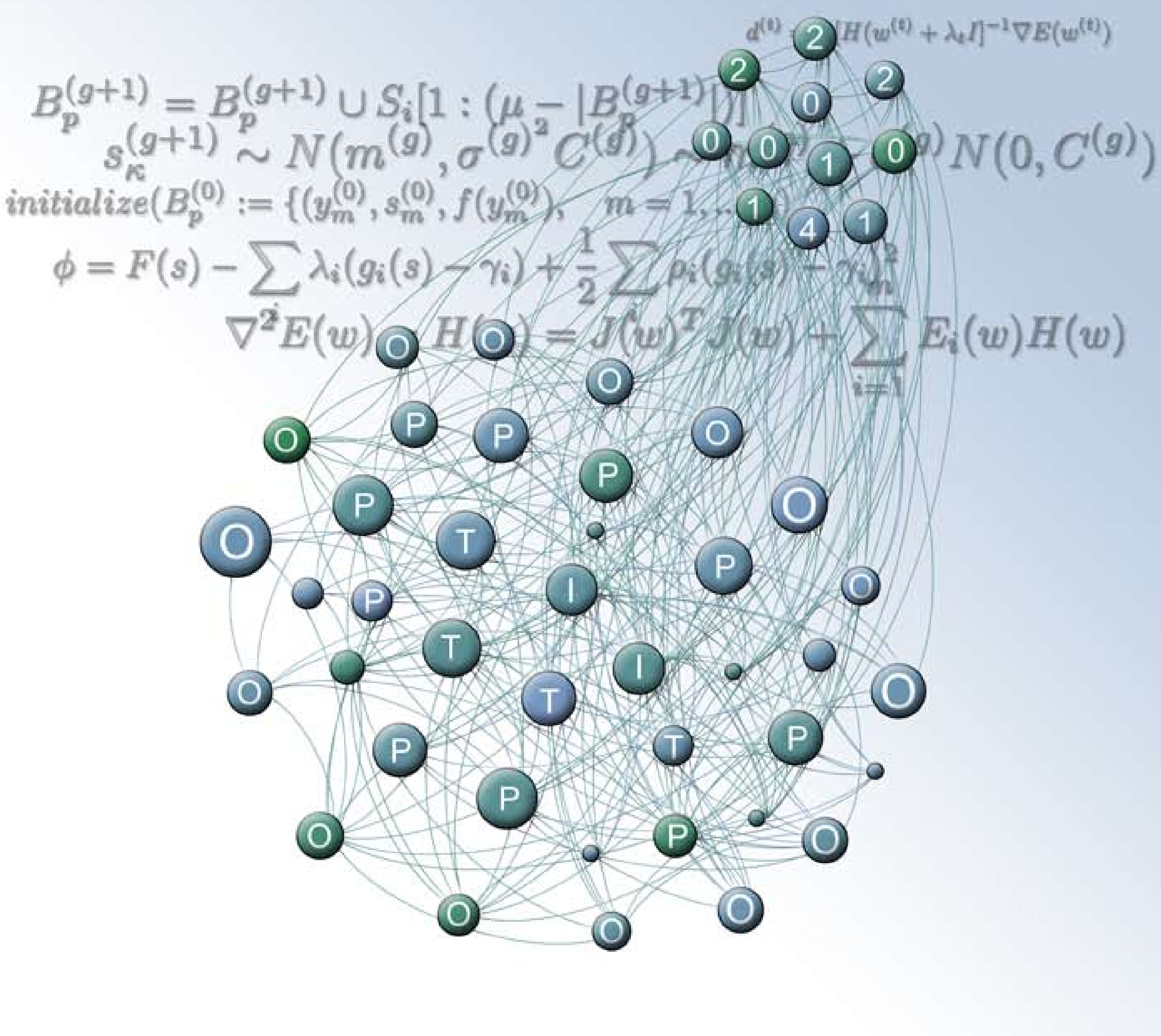

`

# Optimal PMU Placement Using Nonlinear Programming

Nikolaos P. Theodorakatos[1], Nikolaos M. Manousakis[2], George N. Korres[3]

[1] National Technical University of Athens (NTUA)
Iroon Polytechneiou 9, Zografou 15780, Athens, Greece
nikos.theo2772002@gmail.com
[2] National Technical University of Athens (NTUA)
Iroon Polytechneiou 9, Zografou 15780, Athens, Greece
manousakis_n@yahoo.gr; gkorres@cs.ntua.gr

**Keywords:** Phasor Measurement Unit; Optimal Placement; Binary Integer Programming; Mixed Integer Linear Programming; Non Linear Programming; Sequential Quadratic Programming

**Abstract.** Phasor Measurement Units (PMUs) are essential measuring devices for monitoring, control and protection of power systems. The objective of the optimal PMU placement (OPP) problem is to minimize the number of PMUs and select the bus locations to make a power system completely observable. In this paper, the OPP problem is formulated as a nonlinear programming (NLP) problem and a sequential quadratic programming (SQP) method is used for its solution. Simulations are carried out on IEEE standard test systems, using MATLAB. The numerical results are compared to those obtained by a binary integer programming (BIP) model, also implemented in MATLAB. The comparative study shows that the proposed formulation yields the same number of PMUs as the BIP model. The fundamental contribution of this paper lies in investigating the feasibility of using NLP for the solution of the OPP problem and the ability of the proposed methodology to provide multiple solutions in contrast to the binary integer programming model. The System Observability Redundancy Index is adopted to further rank the multiple solutions.

## 1 Introduction

Up to now, monitoring and control of power systems is conducted through the supervisory control and data acquisition (SCADA) system. The SCADA system collects the real-time measurements from the remote terminal units (RTU) placed in substations of the power system. Conventional RTU measurements include power flows, power injections, as well as voltage and line current magnitudes. The phase angle can not be easily measured due to technical difficulties associated with the synchronization of measurements at RTUs. Global Positioning System (GPS) helped to overcome these difficulties and led to the development of Phasor Measurement Units.

A PMU equipped with a GPS receiver provides direct measurement of phase angle with respect to a common reference phase angle [1]. The PMU is placed at a bus to observe the voltage phasor at that bus as well as the current phasors through some or all incident lines. The real time data, provided by PMUs, are transmitted over fast communication links and gathered to higher level devices, known as Phasor Data Concentrators (PDCs) [2], whereas the PMU placement at every substation provides direct measurements of the power system states.

However, it is impossible to install a PMU at every bus of the power system due to the high cost of the PMUs and the lack of communication facilities in some substations. Moreover, as a consequence of Ohm's Law, when a PMU is placed at a bus, the adjacent buses are also observed. This implies that a system can be made observable with a smaller number of PMUs than the number of buses. The optimal PMU placement (OPP) problem involves the determination of the minimum number of PMUs and their corresponding locations in order to achieve complete system observability.

In recent years, there has been significant research activity on the OPP problem. The development and utilization of PMUs were first reported in [3] and [4]. An algorithm for finding the minimum number of PMUs, using a simulated annealing (SA) method and graph theory, is developed in [5]. Reference [5] also reports that the minimum number of PMUs, ensuring full observability of a power system, is $1/5$ to $1/3$ of the system buses. A simple nondominated sorting genetic algorithm that finds the best tradeoffs between competing objectives is proposed in [7]. Four different spanning tree methods based on $N$ and $N-1$ security criteria are suggested in [5]. A graph theoretic PMU placement approach for placing PMUs, based on incomplete observability, is presented in [8].

In addition, several discrete optimization techniques, mathematical or heuristic, have been proposed in literature [9]. Integer linear programming (ILP) is the dominant discrete optimization technique used for solving the OPP problem and many studies concerning this issue have been published [10] - [18]. The ILP technique was initially adopted for the OPP problem solution in [10], [11]. Non linear integer programming and topology transformation of the system are applied to get the OPP solution by considering zero injection buses. Integer Quadratic Programming (IQP) [19], Binary Search Algorithm (BSA) [20], Binary Particle Swarm Optimization (BPSO) [21] and Tabu Search Algorithm (TSA) [22], [23] are some other techniques that have also been implemented for solving the OPP. An iterative weighted least squares algorithm with real decision variables to solve the OPP problem, considering solely PMUs, is introduced in [24]. A global optimization algorithm, Tabu search, is proposed to solve the OPP in [24]. In this paper, a nonlinear programming technique is developed to solve the OPP problem following the formulation [24]. A quadratic objective function is minimized subject to equality nonlinear bus constraints, where the decision variables are defined on the bounded set [0,1]. The quadratic function represents the total PMU installation cost, whereas the nonlinear constraints express the network observability conditions.

The main contribution of this paper lies in investigating the feasibility of using NLP for the OPP problem, despite the fact that this problem is discrete in nature. Hence, we develop a binary integer programming model that guarantees convergence to the optimum solution using existing optimization software. The BIP model is used as a comparative reference to demonstrate the efficiency and accuracy of the proposed model.

The remaining sections of the paper are outlined as follows. Section 2 describes the ILP formulation [10] and the proposed NLP-based framework for solving the OPP problem. The implementation details for each optimization model are presented in Section 3. The power systems used for testing the placement models are described in Section 4. Section 5 provides the simulation results and Section 6 concludes the paper.

## 2. PMU placement problem formulation

A PMU placed at a given bus is capable of measuring the voltage phasor of the bus as well as the phasor currents for all lines incident to that bus. Thus, the entire system can be made observable by placing PMUs at strategic buses in the system [10]. The objective of PMU placement is to minimize the number of PMUs in order to achieve full network observability. In fact, the set of buses where the PMUs have to be installed correspond to a dominating set of the graph [13]. A dominating set (or an externally stable set) in a graph G is a set of vertices that dominates every vertex u in G in the following sense : Either u is included in the dominating set or is adjacent to one or more vertices in the dominating set [32]. Hence, minimum OPP problem maps to smallest dominating set problem on the graph [13].

It is assumed that the PMU has enough channels to measure the voltage phasor at the associated bus and the current phasors of all the lines emanating from that bus [10]. Consequently, the voltage phasors of all adjacent buses will be solvable using the monitored phasor currents along the lines incident to that bus and the known line parameters [19]. In this paper, an ILP [10] and a NLP-based formulation are used to get the OPP problem solution.

### 2.1 Integer Programming: Problem Formulation

For an n-bus system, the OPP problem can be formulated as follows [10]:

$$\min_{x} \; J(x) = \sum_{i=1}^{n} w_i \cdot x_i \tag{1}$$

$$s.t. \;\; f(x) = A \cdot x \geq \hat{1} \tag{2}$$

where $x$ is a binary decision variable vector whose the $ith$ entry, $x_i$, is equal to 1 if a PMU is installed at bus $i$; 0 otherwise, $w_i$ is the cost of PMU installed at bus $i$ and $f(x)$ is a vector function, whose entries are non-zero if the corresponding bus voltage is solvable using the given PMU placement set and zero otherwise. The entries of binary connectivity matrix $A$ are defined as:

$$A_{k,m} = \begin{cases} 1, & \text{if } k = m, \text{ or } k \text{ and } m \text{ are connected} \\ 0, & \text{otherwise} \end{cases} \tag{3}$$

whereas $\hat{1}$, is a vector whose entries are all equal to one.

The IEEE -14 bus system shown in Fig.1 is used to illustrate the ILP approach for the PMU placement problem.

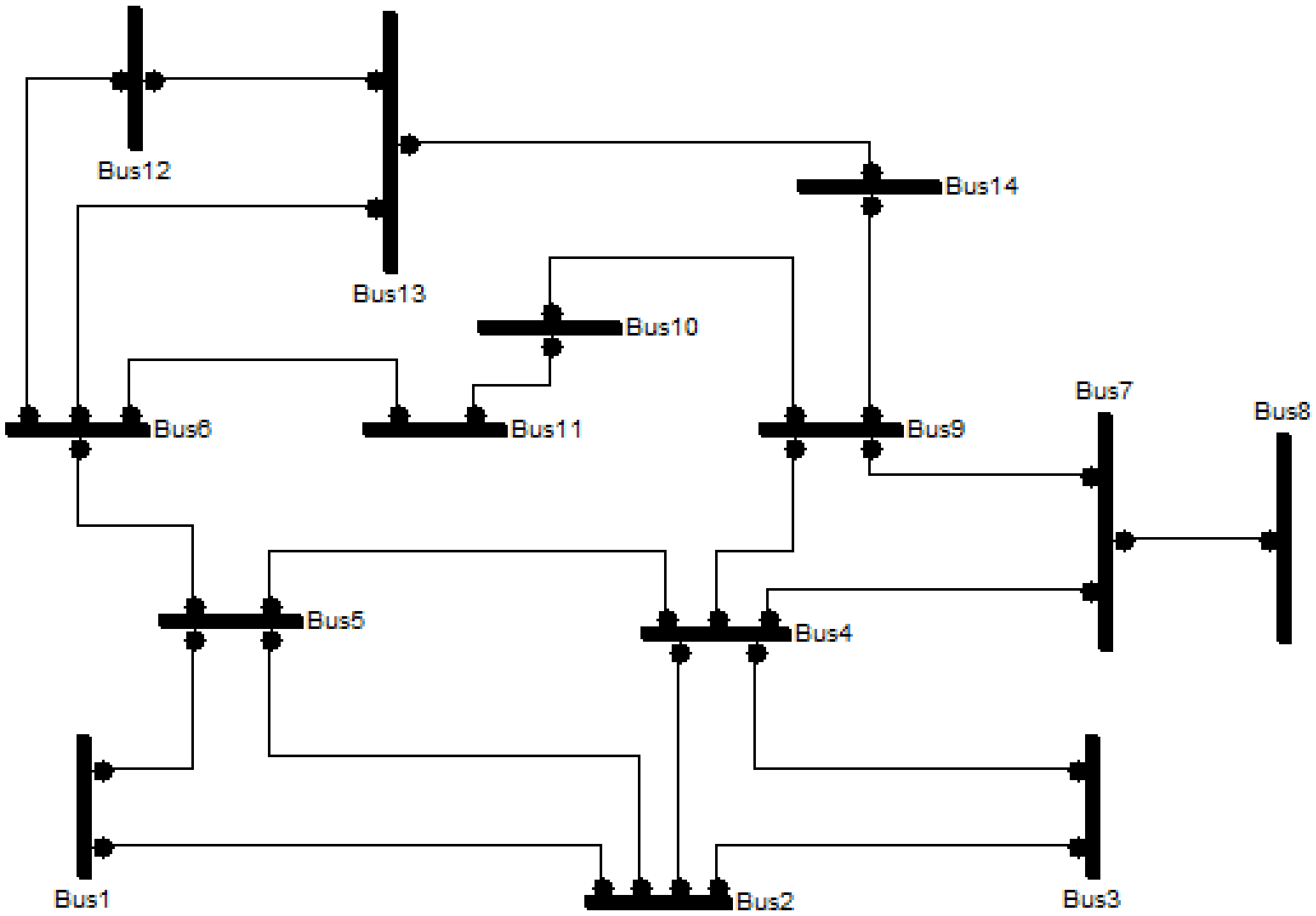

**Figure 1** IEEE 14 -bus system.

The problem formulation is as follows [11]:

$$\min_{x} \; J(x) = \sum_{i=1}^{14} w_i \cdot x_i \tag{4}$$

$$s.t. f(x) = A \cdot x = \begin{cases} f_1 = x_1 + x_2 + x_5 \geq 1 \\ f_2 = x_1 + x_2 + x_3 + x_4 + x_5 \geq 1 \\ f_3 = x_2 + x_3 + x_4 \geq 1 \\ f_4 = x_2 + x_3 + x_4 + x_5 + x_7 + x_9 \geq 1 \\ f_5 = x_1 + x_2 + x_4 + x_5 + x_6 \geq 1 \\ f_6 = x_5 + x_6 + x_{11} + x_{12} + x_{13} \geq 1 \\ f_7 = x_4 + x_7 + x_8 + x_9 \geq 1 \\ f_8 = x_7 + x_8 \geq 1 \\ f_9 = x_4 + x_7 + x_9 + x_{10} + x_{14} \geq 1 \\ f_{10} = x_9 + x_{10} + x_{11} \geq 1 \\ f_{11} = x_6 + x_{10} + x_{11} \geq 1 \\ f_{12} = x_6 + x_{12} + x_{13} \geq 1 \\ f_{13} = x_6 + x_{12} + x_{13} + x_{14} \geq 1 \\ f_{14} = x_9 + x_{13} + x_{14} \geq 1 \end{cases} \tag{5}$$

where: $$x_i \in \{0,\ 1\},\ i = 1, \ldots, 14 \quad (6)$$

### *2.2* Nonlinear Programming: Problem Formulation

Let the continuous decision variable $x_i$ denotes the presence $(x_i = 1)$ or absence $(x_i = 0)$ of a PMU at bus $i$. The OPP problem is formulated as a nonlinear programming problem:

$$\min_x \ J(x) = x^T \cdot W \cdot x = \sum_{i=1}^{n} w_i \cdot x_i^2 \quad (7)$$

$$s.t. \begin{cases} f(x) = 0 \\ \hat{0} \le x \le \hat{1} \end{cases} \quad (8)$$

where $x = (x_1, \cdots, x_n)^T$ is the vector of the decision variables, $J(x) : R^n \rightarrow R$ is the objective function, and $f : R^n \rightarrow R^n$, are the equality observability constraints. $\hat{0}$ and $\hat{1}$ are vectors whose entries are all zeros and ones, respectively. The objective function $J(x)$ is written in matrix notation as $x^T \cdot W \cdot x$, where the matrix $W \in R^{n \times n}$ is a diagonal weight matrix. The diagonal entries $w_i$ of the weight matrix allow the representation of varying installation cost of the PMUs at different buses. In the general case, the PMU installation cost at all buses is the same, $W = I$, where, $I \in R^{n \times n}$ is the identity matrix. Thus, the minimization of $J(x)$ is equivalent to minimizing the total number of PMUs in the system. $f(x)$ is a vector function whose $ith$ entry defines the observability nonlinear equality constraint for the $ith$ bus:

$$f_i(x) = (1 - x_i) \cdot \prod_{j \in a(i)} (1 - x_j) = 0, \quad \forall i \in \mathfrak{I} \quad (9)$$

where $\mathfrak{I}$ is the set of buses and $a(i)$ is the set of buses adjacent to bus $i$. Each inequality constraint (9) implies that at least one PMU should be installed at any one of the buses $i$ and $j \in a(i)$ to make bus $i$ observable.

The binary (boolean) decision variables of the IP approach [10] are transformed into continuous variables by adding the nonlinear observability equality constraints (9). In this way, a consistent system of equations is formulated whose solution is feasible with respect to each equality constraint (9). Mathematically, the formulation (7)−(9) poses no problems to converge to a local optimal solution since all components of $f(x)$ are twice-continuously differentiable. The optimal values of decision variables $x_i$ will be either 1 or 0, as can be proven in Appendix.

The feasible set $S = \{x \mid f_i(x) = 0, 0 \le x_i \le 1,\ i = 1 \ldots n\}$ of the problem is nonconvex. This is because it is made up from equality constraints ( $f_i(x) = 0$ ) which are nonlinear [29]. Because of this, the proposed model is non-convex and can give multiple solutions having the same number of PMUs to the OPP problem solving which they are local minimizers of the optimization problem (7)-(9). Therefore, the optimization problem (7)-(9) have a number of distinct local minimizers. An effective way to obtain these local minimizers in this problem, is to tackle the problem by using a sequential quadratic programming algorithm [26]-[29]. To illustrate the proposed formulation, we use again the IEEE 14- bus system. The OPP problem is formed as follows:

$$\min_{x} J(x) = \sum_{i=1}^{14} x_i^2 \qquad (10)$$

subject to the bus observability constraints:

$$f(x) = \begin{cases} f_1 = (1-x_1)(1-x_2)(1-x_5) = 0 \\ f_2 = (1-x_2)(1-x_1)(1-x_3)(1-x_4)(1-x_5) = 0 \\ f_3 = (1-x_3)(1-x_2)(1-x_4) = 0 \\ f_4 = (1-x_4)(1-x_2)(1-x_3)(1-x_5)(1-x_7)(1-x_9) = 0 \\ f_5 = (1-x_5)(1-x_1)(1-x_2)(1-x_4)(1-x_6) = 0 \\ f_6 = (1-x_6)(1-x_5)(1-x_{11})(1-x_{12})(1-x_{13}) = 0 \\ f_7 = (1-x_7)(1-x_4)(1-x_8)(1-x_9) = 0 \\ f_8 = (1-x_8)(1-x_7) = 0 \\ f_9 = (1-x_9)(1-x_4)(1-x_7)(1-x_{10})(1-x_{14}) = 0 \\ f_{10} = (1-x_{10})(1-x_9)(1-x_{11}) = 0 \\ f_{11} = (1-x_{11})(1-x_6)(1-x_{10}) = 0 \\ f_{12} = (1-x_{12})(1-x_6)(1-x_{13}) = 0 \\ f_{13} = (1-x_{13})(1-x_6)(1-x_{12})(1-x_{14}) = 0 \\ f_{14} = (1-x_{14})(1-x_9)(1-x_{13}) = 0 \end{cases} \qquad (11)$$

where: $$0 \le x_i \le 1,\ i = 1...14 \qquad (12)$$

## 3 Development of PMU placement methodologies

### 3.1 Development of BIP model

The main elements in the BIP model are

1. **Data**.
   - $\Im$ : the set of buses.
   - $n$ : the number of buses.
   - $w_i$ : the weight of the bus $i$
   - $a(i)$ : the set of buses connected through lines to bus $i$.
2. **Variables**. The decision variables involved in this problem are
$$x_i = \begin{cases} 1 & \text{if a PMU is installed at bus } i \\ 0 & \text{otherwise} \end{cases} \qquad (13)$$
3. **Constraints.** The observability inequality constraints are
$$\sum_{j \in a(i)} a_{ij} x_i \ge 1,\ \forall i \in \Im \qquad (14)$$
4. **Function to be minimized.** The total cost is
$$J(x) = \sum_{i=1}^{n} w_i x_i \qquad (15)$$
subject to constraints (13) - (14).

Two solution techniques can be used to solve the binary integer programming model (13)-(15): *the branch-and-bound* (BB) and *branch-and-cut* (BC) methods [25]- [26], [30] - [31]. The BB is the most frequently used and usually the most computationally efficient solution technique [31]. The implementation of BB, provided by the *bintprog* routine of MATLAB, is used to run the BIP model [33]. Fig.2 depicts the BIP flowchart. The optimization problem is solved through the following steps:

Step 1: Read the network branch/bus data.
Step 2: Form the binary connectivity matrix and the PMU cost coefficient vector.
Step 3: Form the right- hand side unity vector.
Step 4: Solve the BIP problem.

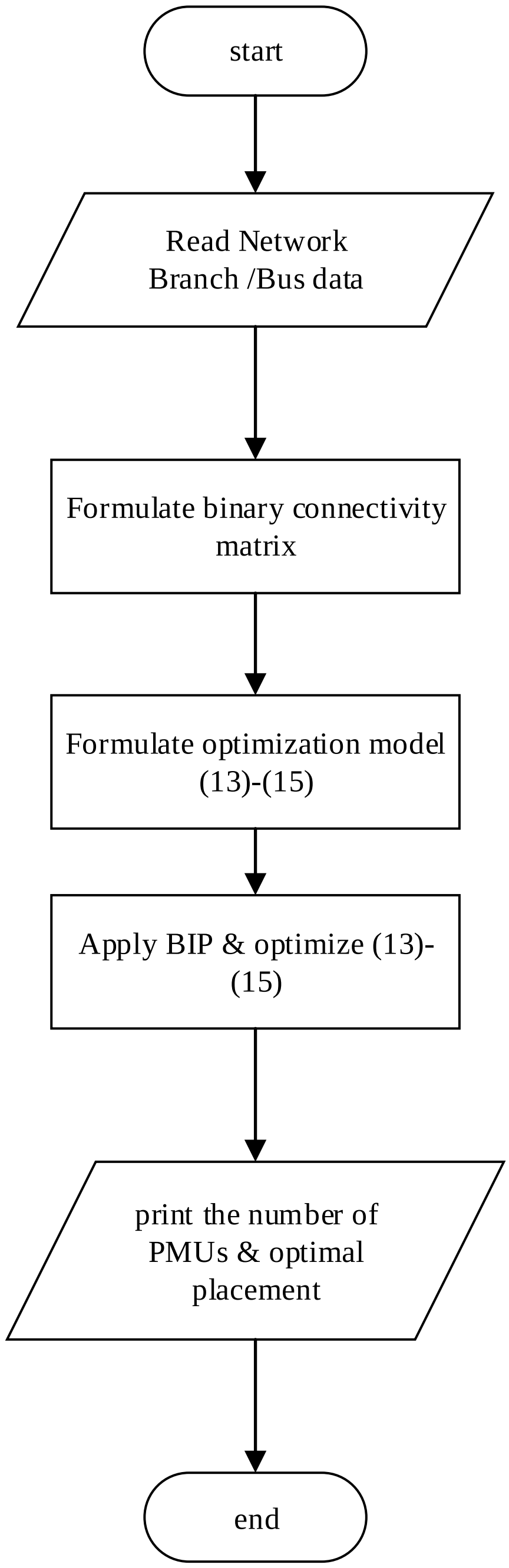

**Figure 2** A flowchart of the BIP-based method for solving the OPP problem

However, an efficient technique of BB, denominated *branch-and-cut* (BC) [31], can be applied to obtain the OPP problem solution with the BIP formulation. Hence, a mixed integer linear programming (MILP) solver named CBC of OPTI Toolbox, an optimization library compatible with MATLAB, can be used to run the BIP model. This solver uses a branch-and-cut algorithm for solving the BIP program [34]. The OPP is formulated as a MILP formulation. The decision variables in the MILP formulation can take integer values [31]. To specify the binary $(0,1)$ variable, first the decision variables $x_i$ are defined to be integer. Then, two constraints are added to specify that the decision variables must be nonnegative and less than or equal to 1. Consequently, the denominated 0/1 MILP formulation is:

$$\min_{x} \ J(x) = \sum_{i=1}^{n} w_i \cdot x_i \tag{16}$$

$$s.t. \begin{cases} A \cdot x \geq \hat{1} \\ x_l \leq x \leq x_u \end{cases} \tag{17}$$

where $x_l$, $x_u$ are the $(n \times 1)$ lower and upper bounds defined as:

$$x_l = [0 \ 0 \ldots.0]^T \tag{18}$$

$$x_u = [1 \ 1 \ldots.1]^T \tag{19}$$

## 3.2 Development of NLP model

The main elements of the NLP model are

1 **Data.**

$\Im$ : the set of buses.

$n$ : the number of buses.

$w_i$ : the weight of the bus $i$

$a(i)$: the set of buses connected through lines to bus $i$.

2. **Variables.** The decision variable vector $x$ is defined on the bounded set.

$$x_l \leq x \leq x_u, \forall i \in \Im \tag{20}$$

where $x_l$, $x_u$ are the $(n \times 1)$ low and upper decision variable bounds defined as:

$$x_l = [0 \ 0 \ldots.0]^T \tag{21}$$

$$x_u = [1 \ 1 \ldots.1]^T \tag{22}$$

3. **Constraints.** The observability equality constraints are

$$f_i(x) = (1 - x_i) \cdot \prod_{j \in a(i)} (1 - x_j) = 0, \ \forall i \in \Im \tag{23}$$

4. **Function to be minimized**. The total cost is

$$J(x) = \sum_{i=1}^{n} w_i \cdot x_i^2 \tag{24}$$

subject to constraints (20) and (23).

The nonlinear programming model (20)-(24) is solved with the *fmincon* NLP solver of MATLAB optimization toolbox. This solver uses a sequential quadratic programming algorithm to solve the constrained minimization problem. We write two m-files to implement the NLP problem in MATLAB [33]. To invoke the objective function by the fmincon, an m-file that returns the current value of the function is written. Another m-file returns the value at the observability constraints at the current $x$ . The decision variables vector $x$ is restricted within certain limits by specifying simple bound constraints to the constrained optimizer routine. The *fmincon* is then executed with a given starting point. The flowchart of the NLP program is shown in Fig.3. The optimization problem is solved through the following steps:

Step 1: Form the objective function.

Step 2: Read the network topology and print the observability constraints.

Step 3: Choose a starting point in the iterative process.

Step 4: Solve the NLP problem.

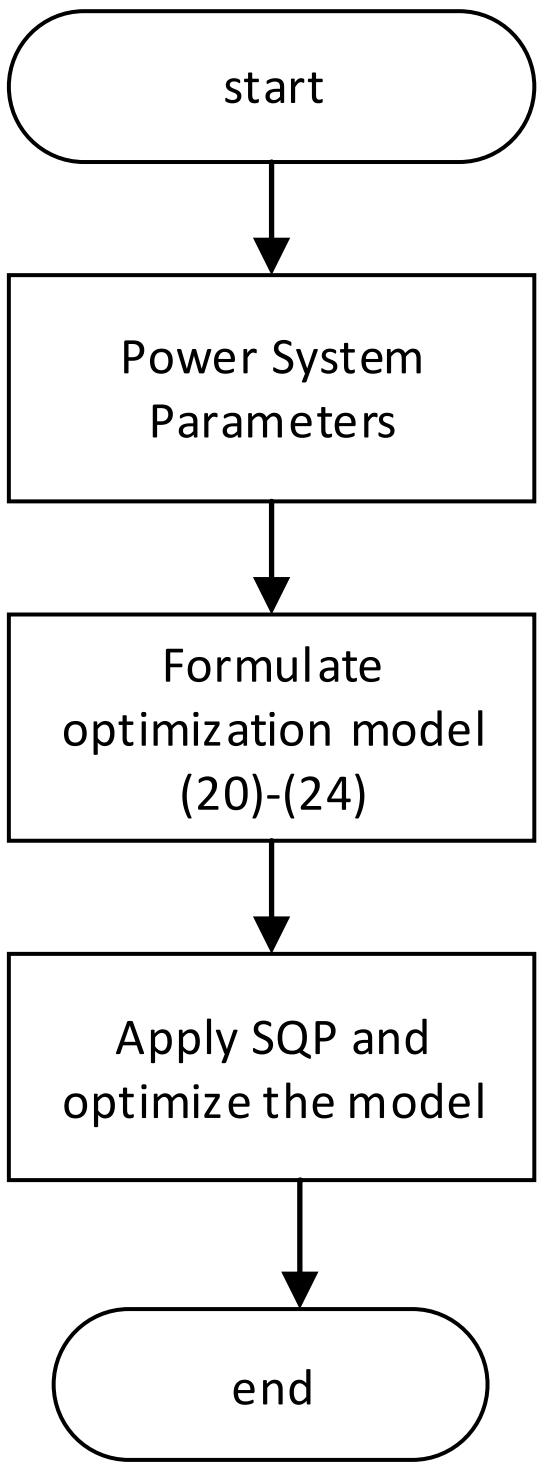

Fig. 3 Flowchart of proposed method for solving the OPP problem.

## 4 Case studies

Power systems differ in terms of the number of buses and the network topology and this makes the task to examine the suitability of a placement methodology with respect to the network size and topology crucial. The developed PMU placement methodologies require the same information, lists of buses and branches, in roughly the same format. For comparison purposes, the PMU placement models are applied to IEEE standard test systems [35]. The characteristics of these test systems are given in Table 1. The numbering of the buses in the IEEE 300 is not successive. The buses of the power system must be re-numbered from 1 up to the total number of buses before the simulation run of each optimization model.

| Test system | Reference | No of buses | No of branches |
|---|---|---|---|
| IEEE 14 bus system | [35] | 14 | 20 |
| IEEE 30 bus system | [35] | 30 | 41 |
| IEEE 57 bus system | [35] | 57 | 80 |
| IEEE 118 bus system | [35] | 118 | 186 |
| IEEE 300 bus system | [35] | 300 | 411 |

**Table 1** General characteristics of the test systems.

## 5 Simulation results and discussion

The computations were carried out using MATLAB optimization solvers. Table 2 summarizes the MATLAB optimization solvers characteristics being used for solving the OPP problem.

| | ILP [10] | | NLP |
|---|---|---|---|
| **Problem Formulation** | **BILP [33]** | **MILP [34]** | |
| **Nature of decision Variables** | Discrete | | Continuous |
| **Decision search** space | $x_i \in \{0,1\}$ | $0 \le x_i \le 1$ | $0 \le x_i \le 1$ |
| **Programming Environment** | MATLAB | | |
| **Solver/Algorithm** | bintprog/ BB | CBC/ branch-and-cut | fmincon/ SQP |
| **Programming solution technique** | ***LP-relaxation*** problem where the binary integer requirement on the variables is replaced by the weaker constraint $0 \le x \le 1$ | Branch and cut involves running a branch and bound algorithm and using cutting planes to tighten the linear programming (LP) relaxations. | SQP methods solve a sequence of optimization subproblems, each of which optimizes a quadratic model of the objective subject to a linearization of the constraints. |

**Table 2** Optimization models used to the OPP problem solving. The following abbreviations are used: ILP = integer linear programming; BB = branch-and-bound; SQP = sequential quadratic programming; LP = linear programming; QP = quadratic programming.

The performance of the proposed model is assessed with respect to the computational time and network size, as well as its ability to consistently provide an acceptable optimum. The NLP optimizer tolerances TolX, TolFun, and TolCon are set, by default, equal to $10^{-6}$ [33], whereas the initial values of the decision variables are set equal to 1, $x_i^0 = 1,\ \forall i \in \Im$ . Furthermore, we set the lower and upper bounds of the decision variables in the NLP solver. The placement results delivered by the *fmincon* are compared with those obtained by using the BIP model, in terms of finding minimum number of PMUs and speed of convergence. To solve the BIP model, *bintprog* requires a feasible point to start. If the starting point is not binary integer feasible, the solver uses the default initial point [33]. The PMU installation weights of each placement model are set equal to 1, $w_i = 1, \forall i \in \Im$ .

The simulation results for the OPP problem are summarized in Table 3. The “Best value” columns present the objective value of the best solution obtained by each optimization solver. From the results, it is obvious that both placement models yield the same minimum number of PMUs and ensure the systems observability. The performance results reveal that, on average, BIP solver employs 0.2956s while MILP solver consumes only 0.0986 s. On the other hand, the NLP solver requires more computational time in comparison to the other solvers, to reach the optimal solution. The computational time, however, is not a serious issue since the PMU placement is a planning problem in nature as it is pointed out in [19].

| Test System | Best Value | BILP | MILP | NLP | |
|---|---|---|---|---|---|
| | | bintprog | CBC | fmincon | |
| | | CPU time (s) | CPU time (s) | Best Value | CPU time (s) |
| IEEE 14 bus | 4 | 0.007 | 0.010 | 4 | 0.060 |
| IEEE 30 bus | 10 | 0.016 | 0.007 | 10 | 0.110 |
| IEEE 57 bus | 17 | 0.155 | 0.020 | 17 | 0.320 |
| IEEE 118 bus | 32 | 0.136 | 0.010 | 32 | 4.050 |
| IEEE 300 bus | 87 | 1.164 | 0.446 | 87 | 28.185 |
| Average | - | 0.2956 | 0.0986 | - | 6.545 |

**Table 3** Optimal number of PMUs obtained by the proposed NLP and ILP methods and required CPU time.

The optimal PMU locations obtained by the ILP solvers, are provided in Tables 4 and 5, respectively. For the IEEE 30-, 57-, 118- bus systems, although the corresponding number of PMUs found by the BIP model is the same as the optimal one reported in [11], the PMU placement set is different.  It is interesting to note that, for a given test system, the ILP solvers deliver different PMU configurations having the same minimum number.

| Test System | PMU location (Bus #) |
|---|---|
| IEEE-14 bus | 2,6,7,9 |
| IEEE-30bus | 1,7,9,10,12,18,24,25,27,28 |
| IEEE-57 bus | 1,4,6,13,19,22,25,27,29,32,36,39,41,45,47,51,54 |
| IEEE-118bus | 3,7,9,11,12,17,21,25,28,34,37,41,45,49,53,56,62,63,68,70,71,76,79,85,86,89,92,96,100, 105,110,114 |
| IEEE-300 bus | 1,2,3,11,12,15,17,22,23,25,26,27,33,37,38,43,48,49,53,54,55,58,59,60,62,64,65,68,71,73, 79,83,85,86,88,92,93,98,99,101,109,111,112,113,116,118,119,128,132,135,138,139,143, 145,152,157,163,167,173,183,187,188,189,190,193,196,202,204,208,210,211,213,216,217, 219,222,226,228,267,268,269,270,272,273,274,276,294 |

**Table 4** Optimal PMU locations obtained by using the bintprog solver

| Test System | PMU location (Bus #) |
|---|---|
| IEEE-14 bus | 2,6,7,9 |
| IEEE-30bus | 1, 2, 6, 9, 10, 12, 18, 24, 25, 27 |
| IEEE-57 bus | 2, 6, 12, 19, 22, 25, 26, 29, 32, 36, 38, 39, 41, 45, 46, 50, 54 |
| IEEE-118 bus | 1, 5, 9, 12, 15, 17, 20, 23, 28, 30, 35, 40, 43, 45, 49, 52, 56, 62, 64, 68, 71, 75, 77, 80, 85, 86, 90, 94, 101, 105, 110, 114 |
| IEEE-300 bus | 1, 2, 3, 11, 12, 13, 15, 17, 22, 23, 25, 27, 29, 33, 37, 38, 41, 43, 48, 49, 53, 54, 55, 58, 59 60, 62, 64, 65, 68, 71, 76, 83, 85, 86, 88, 93, 98, 99, 101, 103, 109, 111, 112, 113, 116 118, 119, 122, 132, 135, 138, 143, 145, 152, 157, 163, 167, 168, 173 183, 187, 189, 190 193, 196, 200, 204, 208, 210, 211, 213, 216, 217, 219, 222, 225, 228, 267, 268, 269, 270 272, 273 274, 276, 294 |

**Table 5** Optimal PMU locations obtained by using the CBC solver.

Table 6 provides the optimal PMU locations obtained by the NLP solver. The PMU placement sets are different from those found by using the BIP model. These results confirm the observation reported in [11],[19] that there can be more than one solution to the OPP problem with the same cost.

| Test System | PMU location (Bus #) |
|---|---|
| IEEE-14 bus | 2,8,10,13 |
| IEEE-30bus | 1, 2, 6, 9,10,12,15,20,25,27 |
| IEEE-57 bus | 1, 4, 9,15 ,20 ,24 ,25,28 ,29 ,32 ,36 ,38,41, 46 ,50,53, 57 |
| IEEE-118 bus | 2, 5, 9, 12, 15, 17, 21, 23, 25, 28, 34, 37, 40, 45, 49, 52, 56, 62, 64, 68, 71, 75, 77, 80,85, 87, 91, 94, 101, 105, 110,114 |
| IEEE-300 bus | 1, 2, 3, 11, 12, 15, 17, 22, 23, 25, 26, 27, 29, 33, 37, 38, 43, 48, 49, 53, 54, 59, 62, 64,65, 68, 71, 79, 82, 85, 86, 88, 89, 93, 98, 99, 101, 109, 111, 112, 113, 116, 118, 119, 124, 132, 135, 138, 139, 143, 145, 152, 157, 163, 167, 173, 177, 183, 187, 189, 190, 193, 196, 202, 204, 209, 210, 212, 213, 216, 217, 224, 225, 228, 230, 236, 237 ,238, 267, 268, 269, 270, 272, 273, 274, 276, 294 |

**Table 6** Optimal PMU locations obtained by using the fmincon solver.

Another issue investigated in this paper is the starting point selection for solving the described optimization models. Starting from the default initial point, the ILP solvers can get only one optimal placement set, whereas more than one solution may exist [11], [19]-[20]. Instead, the NLP model may yield more than one optimal solutions with the same minimum number of PMUs. To get more than one optimal solutions, we solve the proposed model with the *fmincon*, starting from different initial points selected within the variable bounds $[0,1]$. The selection of a starting point can be made for example with a step of 0.1 among the selected initial points. Therefore, any point which belongs to the feasible set ( $x \in S$ ) can be chosen for an initial design starting point, where the feasible set is $S = \{ x \mid f_i(x) = 0, 0 \le x_i \le 1, \ i = 1...n \}$. The only difference between two starting points is that the selection of the initial point may affect the number of iterations in order to converge to a local minimum of the minimization problem [25]- [29]. In any case, we have found that a choice of any starting point leads to a distinct local minimizer of the NLP problem. We have also shown that the NLP results agree with those found by the BILP model regarding the number of PMUs required for full *system* observability (Table 3). Multiple solutions exist for the test cases shown in Table 7. We adopt the system observability redundancy index (SORI) [13], to further rank these multiple placement solutions. The solution that maximizes the SORI index is denoted by bold characters.

| Test System | PMU | PMU location (Bus #) | SORI |
|---|---|---|---|
| IEEE 14 bus | 4 | 2,8,10,13 | 14 |
| | | 2 ,6 ,8 ,9 | 17 |
| | | 2,7,11,13 | 16 |
| | | 2,7,10,13 | 16 |
| | | **2,6,7,9** | **19** |
| IEEE 30 bus | 10 | **1, 2, 6, 9,10,12,15,20,25,27** | **50** |
| | | 1, 5,9,10,12,15,18,25,28,29 | 42 |
| | | 2,4,6,9,10,12,18,24,26,29 | 47 |
| | | 2,3,6,9,10,12,18,24,25,29 | 47 |
| | | 3,5,9,10,12,15,19,25,27,28 | 44 |
| | | 2,4,6,10,11,12,15,18,25,29 | 48 |
| | | 3,5,8,9,10,12,18,23,26,30 | 37 |
| | | 3, 5, 6 ,9 ,10,12,15,20,25,29 | 46 |
| | | 1, 5, 8,10,11,12,19,23,26,27 | 37 |
| | | 1,5, 8 ,10,11,12,19 ,23,26 ,29 | 35 |
| | | 3,5,8,10,11,12,18,23,26,29 | 35 |
| | | 3, 5, 8, 9 ,10,12,15,18,25,30 | 41 |

**Table 7** Multiple optimal solutions obtained by using the fmincon solver, for different starting points.

| Test System | PMU | PMU location (Bus #) | SORI |
|---|---|---|---|
| IEEE 30 bus | 10 | 3,5,8,9,10,12,19,23,26,30 | 37 |
| | | 1,7 ,8 ,10 ,11,12,19,23,26 ,29 | 35 |
| | | 1,7,9,10,12,15,20,25,28,30 | 42 |
| | | 3,5,8,10,11,12,18,23,25,29 | 37 |
| | | 1,7,8,9,10,12,15,19,25,29 | 41 |
| | | 3,5,8,10,11,12,18, 24, 25, 30 | 38 |
| | | 1,5,9,10,12,18,23,26,28,30 | 38 |
| | | 1,5,8,10,11,12,18,24,26,29 | 36 |
| IEEE 57 bus | 17 | 1,4,9,15,20,23,25,27 ,29 ,32,36,41,44,47,50,54,57 | 67 |
| | | 2,6,12,19,22,25,27,29,32,36,41,45,46,47,50,54,57 | 62 |
| | | 2,6,12,19,22,25,27,32,36,39,41,45,46,49,50,52,55 | 63 |
| | | 1, 4, 6,10, 19,22,25, 27,32,36,41,45,46,49,52,55,57 | 65 |
| | | 2,6, 12,19, 22, 27,32, 36, 39,41, 45, 46,47,50,52,55,57 | 61 |
| | | **1, 4, 9,15 ,20 ,24 ,25,28 ,29 ,32 ,36 ,38,41, 46 ,50,53, 57** | **71** |
| | | 1,6 ,7,9,15,19,22,25,27,32,36,38,39,41,47,50,53 | 71 |
| | | 1, 4, 7, 9, 15, 19, 22, 25, 27, 32, 36, 38, 41, 47, 50, 53, 57 | 71 |
| | | 1,3,6,10,19,22,25,27,32,36,41,44,46,49,52,55,57 | 64 |
| | | 1,4,6,10,19,22,25,27,32,36,41,45,46,49,52,55,57 | 65 |
| | | 1,6,10,15,19,22,25,27,32,36,41,44,46,49,52,55,57 | 66 |
| | | 1,6,10,15,19,22,25,27,32,36,41,45,46,49,52,55,57 | 66 |
| | | 1,4,6,10,19,22,25,27,29,32,36,41,44,46,49,54,57 | 66 |
| | | 3,6,12,15,19,22,25,27,32,36,38,39,41,46,50,52,55 | 68 |
| | | 1,4,9,10,19,22,25,26,29,32,36,41,44,46,49,54,57 | 68 |
| | | 1,4,9,13,19,22,26,29,30,32,36,39,41,45,47,51,54 | 69 |
| | | 1,6,10,15,19,22,25,27,32,36,38,41,46,49,52,55,57 | 69 |
| | | 2,6,12,19,22,25,27,32,36,41,45,46,47,50,52,55,57 | 61 |
| | | 1,4,9,13,19,22,26,29,30,32,36,39,41,44,47,50,54 | 69 |
| | | 1 ,4,9,20,22,25,27,29,32,36 ,41,44 ,46 ,49 ,50,53,57 | 67 |
| IEEE 118 bus | 32 | 2, 5, 9, 12, 15, 17, 21, 23, 25, 28, 34, 37, 40, 45, 49, 52, 56, 62, 64, 68, 71, 75, 77, 80,85, 87, 91, 94, 101, 105, 110, 114 | 161 |
| | | 2,5,9,12,15,17,21,25,29,34,37,41,45,49,52,56,62,64,68,70,71,75,77,80,85,86 91,94,101,105,110,114 | 161 |
| | | 2, 5, 9, 11, 12, 17, 21, 23, 25, 29, 34, 37, 41, 45, 49, 52, 56, 62, 64, 68, 71, 75,77, 80, 85, 86, 90, 94, 102, 105, 110, 115 | 159 |
| | | 1 ,5 ,9 ,11 ,12 ,17 ,21 ,25 ,29 ,34 ,37 ,40 ,45 ,49 ,52 ,56 ,62 ,64 ,68 ,71 ,72 75 ,77 ,80 85, 86 ,91 ,94 ,101 ,105 ,110 ,114 | 159 |
| | | 1 ,5 ,9 ,11 ,12 ,17 ,21 ,25 ,29 ,34 ,37 ,40 ,45 ,49 ,52 ,56 ,62 ,64 ,72 ,73 ,75 77 ,80 ,85,86 ,91 ,94 ,101 ,105 ,110 ,114 ,116 | 154 |
| | | 1 ,5 ,9 ,11 ,12 ,17 ,21 ,25 ,29 ,34 ,37 ,40 ,45 ,49 ,52 ,56 ,62 ,64 ,72 ,73 ,75 77 ,80 ,85,86 ,91 ,94 ,101 ,105 ,110 ,114 ,116 | 152 |
| | | 1 ,5 ,10 ,11 ,12 ,17 ,21 ,25 ,29 ,34 ,37 ,41 ,45 ,49 ,52 ,56 ,62 ,64 ,72 ,73 ,75 77,80, 85, 87 ,91 ,94 ,101 ,105 ,110 ,114 ,116 | 150 |
| | | 1 ,5 ,10 ,12 ,13 ,17 ,21 ,25 ,29 ,34 ,37 ,41 ,45 ,49 ,53 ,56 ,62 ,64 ,72 ,73 ,75 77,80, 85,87 ,91 ,94 ,102 ,105 ,110 ,114 ,116 | 148 |
| | | 2, 5, 9, 12, 15, 17, 21, 25, 29, 34, 37, 40, 45, 49, 52, 56, 62, 64, 68, 71, 72, 75, 77,80, 85, 86, 91, 94, 101, 105, 110, 114 | 160 |
| | | 2, 5, 9, 12, 15, 17, 21, 25, 29, 34, 37, 40, 45, 49, 52, 56, 62, 63, 68, 70, 71, 75, 77,80, 85, 86, 90, 94, 102, 105, 110, 114 | 162 |
| | | 2, 5, 9, 12, 15, 17, 21, 23, 25, 29, 34, 37, 40, 45, 49, 52, 56, 62, 64, 68, 71 75, 77, 80, 85, 86, 91, 94, 101, 105, 110, 114 | 162 |
| | | **2, 5, 9, 12, 15, 17, 21, 25, 29, 34, 37, 40, 45, 49, 52 ,56 , 62 , 64 , 68, 70, 71, 75, 77, 80, 85, 86, 91, 94, 101, 105, 110, 114** | **163** |
| | | 2,5,9,12,15,17,21,25,29,34,37,40,45,49,52,56,62,64,71,72,75,77,80,85,86,90 94,101,105,110,114,116 | 157 |

Table 7(CONTINUED)

| Test System | PMU | PMU location (Bus #) | SORI |
|---|---|---|---|
| IEEE 118 bus | 32 | 1,7,9,11,12,17,21,25,29,34,37,41,45,49,52,56,62,64,72,73,75,77,80,85,87,91 94,101,105,110,114,116 | 148 |
| | | 3,5,9,12,15,17,21,23,25,28,34,37,40,45,49,52,56,62,64,68,71 75,77,80,85,86,91,94,101,105,110,114 | 163 |
| | | 2,5,9,12,15,17,21,23,25,28,34,37,40,45,49,52,56,62,64,71,75,77 80,85,87,91,94,101,105,110,114,116 | 158 |
| | | 3,5,9,11,12,17,21,25,28,34,37,40,45,49,52,56,62,63,68,70,71,75,77,80,85,86 90,94,102,105,110,114 | 162 |
| | | 3,5,9,11,12,17,21,25,29,34,37,40,45,49,52,56,62,63,68,70,71,75,77,80,85,86 90,94,102,105,110,114 | 162 |
| | | 3,5,9,11,12,17,21,23,25,29,34,37,40,45,49,52,56,62,63,68,71,75,77,80,85,86 90,94,102,105,110,114 | 161 |
| IEEE 300 bus | 87 | 1,2,3,11,12,15,17,19,22,23,25,27,33,37,38,41,43,48,49,53,54,62,64,65,68 71,73,79,82,85,86,88,93,98,99,101,109,111,112,113,116,118,119,124,132 135,138,139,141,145,152,157,163,167,173,177,183,187,189,190,193,196 202,204,209,210,212,213,216,217,221,223,228,230,236,237,238,262,267 268,269,270,272,273,274,276,294 | 411 |
| | | 1 ,2 ,3 ,11 ,12 ,15 ,17 ,22 ,23 ,25 ,26 ,27 ,33 ,37 ,38 ,43 ,48 ,49 ,53 54 ,59 ,62 ,64 ,65, 68 ,71 ,73 ,79 , 82 ,85 ,86 , 88 ,89 ,93 ,98 ,99 ,101 103 ,109 ,111 ,112 ,113 ,116 ,118, 119 ,124 ,132 ,135 ,138 ,143 ,145 152 ,157 ,163 ,167 ,173 ,177 ,183 ,187 ,189 ,190,193 ,196 ,202 ,204 209 ,210 ,212 ,213 ,216 ,217 ,221 ,223 ,228 ,230 ,236 ,237 ,238 267 ,268 ,269 ,270 ,272 ,273 ,274 ,276 ,294 | 415 |
| | | 1,2,3,11,15,17,22,23,25,26,27,33,37,38,43,48,49,53,54,59,62,64,65,68,71,73 79,82,85,86,88,89,93,98,99,101,103,109,111,112,113,116,118,119,124,132, 135,138,143,145,152,157,163,167,173,183,187,188,189,190,193,196,202, 204,209,210,212,213,216,217,221,223,228,230,236,237,238,251,267,268, 269,270,272,273,274,276,294 | 408 |
| | | 1,2,3,11,15,17,22,23,25,26,27,33,37,38,43,48,49,53,54,59,62,64,65,68,71,73 79,82,85,86,88,89,93,98,99,101,103,109,111,112,113,116,118,119,124,132, 135,138,143,145,152,157,163,167,173,177,183,187,189,190,193,196,202, 204,209,210,212,213,216,217,221,223,228,230,236,237,238,251,267,268, 269,270,272,273,274,276,294 | 412 |
| | | 1,2,3,11,12,15,17,19,22,23,25,27,33,37,38,43,48,49,53,54,62,64,65,68,71,73 79,82,85,86,88,89,93,98,99,101,109,111,112,113,116,118,119,124,132,135, 138,139,141,145,152,157,163,167,173,177,183,187,189,190,193,196,202, 204,209,210,212,213,216,217,221,223,228,230,236,237,238,262,267,268, 269,270,272,273,274,276,294 | 412 |
| | | 1,2, 3,11,15,17,22,23,25,26,27,29,33,37,38,43,48,49,53,54,59,62 64,65,68,71,79,82,85,86,88,89,93,98,99,101,109,111,112,113,116 118,119,124,132,135,138,139,143,145,152,157,160,163,173 177,183,187,189,190,193,196,202,204,209,210,212,213,216,217 224,225,228,230,236,237,238,251,267,268,269,270,272,273,274,276,294 | 413 |
| | | 1,2,3,11,12,15,17,19,22,23,25,27,33,37,38,43,48,49,53,54,59,62,64,65,68,71 73,79,82,85,86,88,89,93,98,99,101,109,111,112,113,116,118,119,124,132, 135,138,139,143,145,152,157,163,167,173,177,183,187,189,190,193,196, 202,204,209,210,212,213,216,217,221,223,228,230,236,237,238,267,268, 269 270,272,273,274,276,294 | 416 |
| | | 1,2,3,11,12,15,17,22,23,25,26,27,29,33,37,38,43,48,49,53,54,59 62,64,65,68,71,79,82,85,86,88,89,93,98,99,101,109,111,112,113 116,118,119,124,132,135,138,139,143,145,152,157,163,167 173,177,183,187,189,190,193,196,202,204,209,210,212,213 216,217,223,224,228,230,236,237,238,267,268,269,270,272 273,274,276,294 | 416 |

**Table 7(CONTINUED)**

| Test System | PMU | PMU location (Bus #) | SORI |
|---|---|---|---|
| IEEE 300 bus | 87 | 1, 2, 3, 11, 12, 15, 17, 22, 23, 25, 26, 27, 29, 33, 37, 38, 43, 48, 49, 53, 54, 59, 62, 64,65, 68, 71, 79, 82, 85, 86, 88, 89, 93, 98, 99, 101, 109, 111, 112, 113, 116, 118, 119, 124, 132, 135, 138, 139, 143, 145, 152, 157, 163, 167, 173, 177, 183, 187, 189, 190, 193, 196, 202, 204, 209, 210, 212, 213, 216, 217, 224, 225, 228, 230, 236, 237 ,238, 267, 268, 269, 270, 272, 273, 274, 276, 294 | 417 |
| | | **1,2,3,11,12 ,15,17 ,19 ,22, 23,25 ,27, 33,37,38,43,48,49, 53,54,59,62,64,65, 68,71,73, 79,82,85,86,88,89,93,98,99,101,109,111,112, 113,116,118,119,124,132,135,138,139,143,145,152,157,163,167, 173,177,183,187,189,190,193,196,202,204,209,210,212, 213,216,217,224,225,228,230,236,237, 238,267,268,269,270,272,273,274,276,294** | **418** |
| | | 1,2,3,11,15,17,22,23,25,26,27,29,33,37,38,43,48,49,53,54,59,62,64,65,68,71 79,82,85,86,89,93,98,99,101,109,111,112,113,116,119,124,132,135,138,139 143,145,152,157,160,163,173,183,187,188,189,190,193,196,202,204,209 210,212,215,216,217,224,225,228,230,235,236,237,238,251,264,267,268 269,270,272,273,274,276,294 | 402 |
| | | 1,2,3,11,15,17,22,23,25,26,27,29,33,37,38,43,48,49,53,54,59,64,65,68,71,79 82,85,86,88,89,93,98,99,101,109,111,112,113,116,118,119,124,132,135,138 139,143,145,152,157,163,167,173,177,183,187,189,190,193,196,202,204, 209,210,212,213,216,217,222,225,228,230,236,237,238,240,251,267,268 269,270,272,273,274,276,294 | 409 |
| | | 1,2,3,11,15,17,22,23,25,26,27,29,33,37,38,43,48,49,53,54,59,62,64,65,68,71 79,82,85,86,88,89,93,98,99,101,109,111,112,113,116,118,119,124,132,135, 138,139,143,145,152,157,163,167,173,177,183,187,189,190,193,196,202 204,209,210,212,213,216,217,223,224,228,230,236,237,238,251,267,268 269,270,272,273,274,276,294 | 413 |
| | | 1,2,3,11,12,15,17,20,22,23,25,27,29,33,37,38,43,48,49,53,54,59 62,64,65,68,71,79,82,85,86,88,89,93,98,99,101,103,109,111,112 113,116,118,119,124,132,135,138,143,145,152,157,163,167,173 177,183,187,189,190,193,196,202,204,209,210,212,213,216,217 224,225,228,230,236,237,238,267,268,269,270,272,273,274,276,294 | 420 |
| | | 1 ,2 ,3 ,11 ,15 ,17 ,22 ,23 ,25 ,26 ,27 ,33 ,37 ,38 ,41 ,43 ,48 ,49 53 ,54 ,59 ,62 ,64 ,65 ,68 ,71 ,73 ,79 ,82 ,85 ,86 ,88 ,93 ,98 , 99 101 ,103 ,109 ,111 ,112 ,113 ,116 ,118 ,119 ,124 ,132 ,135 ,138 143 ,145 ,152 ,157 ,163 ,167 ,173 ,177 ,183 ,187 ,189 ,190 ,193 196 ,202 ,204 ,209 ,210 ,212 ,213 ,216 ,217 ,221 ,223 ,228 ,230 236 ,237 ,238 ,251 ,267 ,268 ,269 ,270 ,272 ,273 ,274 ,276 ,294 | 411 |

**Table 7(CONTINUED)**

## 6. Conclusions

This paper presents a nonlinear programming model for the OPP problem ensuring the complete system observability. The proposed methodology was implemented in MATLAB, using sequential quadratic programming, and successfully tested on different size power systems. The test results were compared with those obtained by a binary integer programming model implemented in MATLAB, and they validate the effectiveness and accuracy of the NLP model. Depending upon the starting point, the developed optimization scheme is able to yield different PMU placement sets having the same minimum number of PMUs. The proposed PMU placement method ensures the power system observability in the absence of any conventional measurement. Future work will include additional constraints into the proposed model, such as the existence of zero injection, and power flow measurements.

## Appendix

Consider the nonlinear equality constraints $f_i(x) = (1 - x_i) \cdot \prod_{j \in a(i)} (1 - x_j) = 0, \ \forall i \in \Im$. The optimization problem can be stated as follows:

$$\min_{x \in R^n} \{ J(x) : \ f_i(x) = 0, \ 0 \le x_i \le 1, \ i = 1, \ldots, n \} \tag{A.1}$$

Suppose that point $x^*$ is a local minimizer of the optimization problem and there exists a $k \in \{1, \ldots, n\}$ such that:

$$x_k^* \in (0,1) \tag{A.2}$$

In addition, we have that:

$$f_i\left(x^*\right) = 0, \ i = 1, \ldots, n \tag{A.3}$$

$$0 \le x_i^* \le 1, \ i = 1, \ldots, k-1, k+1, \ldots, n \tag{A.4}$$

Equations (A.3) are satisfied at the point $x^*$, when the terms $\left(1 - x_i^*\right), \ i \ne k,$ become equal to zero (the term $\left(1 - x_k^*\right)$ is non-zero). These terms are sufficient to satisfy equations (A.3), $\forall \ i \in \{1, \ldots, n\}$. Hence, the points $\hat{x}(\delta) = x^* + \delta \cdot e_k, \ \forall \ \delta \in R$ also satisfy the equations:

$$f_i\left(\hat{x}(\delta)\right) = 0, \ \forall \ i = \{1, \ldots, n\}, \forall \ \delta \in R \tag{A.5}$$

Moreover, we have that:

$$\left.\begin{aligned} J\left(\hat{x}(\delta)\right) &= \hat{x}(\delta)^T \cdot W \cdot \hat{x}(\delta) = \left(x^* + \delta \cdot e_k\right)^T \cdot W \cdot \left(x^* + \delta \cdot e_k\right) \\ &= x^{*T} \cdot W \cdot x^* + 2 \cdot \delta \cdot e_k^T \cdot W \cdot x^* + \delta^2 \cdot e_k^T \cdot W \cdot e_k \\ &= J\left(x^*\right) + 2 \cdot \delta \cdot w_k \cdot x_k^* + \delta^2 \cdot w_k \\ &= J\left(x^*\right) + \delta \cdot w_k \cdot \left(2x_k^* + \delta\right) \end{aligned}\right\} \Rightarrow J\left(\hat{x}(\delta)\right) < J\left(x^*\right), \ \forall \ \delta \in \left(-2x_k^*, 0\right) \tag{A.6}$$

and

$$\begin{aligned} &0 \le \hat{x}_i(\delta) = x_i^* \le 1, \ \forall \ i \in \{1, \ldots, n\} \setminus \{k\} \\ &0 \le \hat{x}_k(\delta) = x_k^* + \delta \le 1, \ \forall \ \delta \in \left[-x_k^*, 1 - x_k^*\right] \end{aligned} \tag{A.7}$$

From (A.5)–(A.7), we can conclude that the points $\hat{x}(\delta), \ \forall \ \delta \in [-x_k^*, 0)$ satisfy all the constraints of the above optimization problem and that $J\left(\hat{x}(\delta)\right) < J\left(x^*\right), \ \forall \ \delta \in [-x_k^*, 0)$. Therefore, the point $x^*$ is not a local minimum of the optimization problem.

Given a local minimum $x^*$ of (A.1), some of the $x_k^*$ will be equal to $1$ satisfying the observability constraints $f_i(x^*) = 0, \ i = 1, \ldots, n$. The rest of $x_i^*, \ i \ne k,$ will be equal to $0$, because the cost function $\sum_{i=1}^{n} w_i \cdot x_i^2$, given $0 \le x_i \le 1, \ i = 1, \ldots, n,$ is minimized when $x_i^* = 0$.

ACKNOWLEDGMENT

Nikolaos P. Theodorakatos would like to thank his Prof. Nicholas Maratos at School of Electrical and Computer Engineering at NTUA for the teaching of optimization techniques and giving the proof of appendix.

**OPTI-*i***
**Proceedings of the**
**1st International Conference on**
**Engineering and Applied Sciences Optimization**

M.G. Karlaftis, N.D. Lagaros, M. Papadrakakis (Eds.)